\UseRawInputEncoding
\documentclass[12pt, leqno]{article}
\usepackage{latexsym}
\usepackage{enumerate,amsmath,amsthm,amssymb,amsfonts,amscd}

\newtheorem{defn}{Definition}[section]
\newtheorem{thm}[defn]{Theorem}
\newtheorem{cor}[defn]{Corollary}
\newtheorem{prop}[defn]{Proposition}

\newtheorem{rem}[defn]{Remark}

\numberwithin{equation}{section}

\def\A{\mbox{${\cal A}$}}
\def\B{\mbox{${\cal B}$}}
\def\F{\mbox{${\cal F}$}}
\def\U{\mbox{${\cal U}$}}
\def\O{\mbox{${\cal O}$}}

\newcommand{\rr}{\mathbb R}
\newcommand{\cc}{\mathbb C}
\newcommand{\zz}{\mathbb Z}
\newcommand{\nn}{\mathbb N}

\renewcommand{\a}{\alpha}
\newcommand{\be}{\beta}
\newcommand{\ga}{\gamma}
\newcommand{\vp}{\varphi}
\newcommand{\vr}{\varrho}
\newcommand{\vep}{\varepsilon}

\newcommand{\ra}{\rightarrow}
\newcommand{\la}{\leftarrow}
\newcommand{\lra}{\longrightarrow}
\newcommand{\lla}{\longleftarrow}

\input{arrow.tex}

\begin{document}
\title{\bf On avatars of the Singer-Wermer conjecture in Fr\'{e}chet algebras}

\author{S. R. PATEL}
\maketitle
\begin{tabbing}
\noindent {\bf 2010 Mathematics Subject Classification:} \= Primary 46J05;\\
\> Secondary 13F25, 46H40
\end{tabbing}
\noindent {\bf Key words:} Fr\'{e}chet algebra of power series in
infinitely many indeterminates, derivations, (in)equivalent
Fr\'{e}chet algebra topologies, Singer-Wermer conjecture.\\

\newpage
\noindent {\bf Abstract.} We construct two Fr\'{e}chet algebras
admitting countably many mutually inequivalent Fr\'{e}chet algebra
topologies. The second example is a modification of the maiden
(and first) example of a non-Banach Fr\'{e}chet algebra with two
inequivalent Fr\'{e}chet algebra topologies, considered by Read to
show that the famous Singer-Wermer conjecture does not hold in the
Fr\'{e}chet case. We discuss a couple of avatars of the
Singer-Wermer conjecture in Fr\'{e}chet algebras. In particular,
this conjecture (and so, its weaker version, too) does hold by the
second example which also admits countably many mutually
equivalent Fr\'{e}chet algebra topologies. However the stronger
version of this conjecture (i.e., the Kaplansky's conjecture) does not hold by the second example.\\
\newpage
\section{Introduction} Recently, the author posed a question about
the existence of a Fr\'{e}chet algebra (i.e., complete, metrizable
locally $m$-convex algebra) with infinitely many inequivalent
Fr\'{e}chet algebra topologies ~\cite[Question 4]{12}. This
question carries a lot of importance due to the following three
reasons.

First, in automatic continuity theory, we would normally like to
examine when and how the algebraic structure of the algebra $A$
determines the topological structure of $A$, in particular, the
continuity aspect (and more particularly, the uniqueness of the
topology of $A$; see ~\cite{3, 6, 10, 11, 12} for more details).
So it is natural to expect that the non-uniqueness of the topology
would reflect some properties of the algebraic structure of $A$.

In ~\cite{7}, Feldman constructed an example to show the failure
of the Wedderburn Theorem in the Banach algebra case. This algebra
is $\ell_2 \,\oplus \,\cc$ as a vector space with the usual
product in $\ell_2$ and trivial multiplication by the second
summand, and has a norm in which $\ell_2$ is dense. In ~\cite{13},
Read constructed an example (see \S 2 below) to show the failure
of the Singer-Wermer conjecture in the commutative Fr\'{e}chet
algebra case; this conjecture holds in the commutative Banach
algebra case (see ~\cite{14}). Indeed, a sharp observation of
these developments leads to a couple of avatars of the
Singer-Wermer conjecture such as a weaker version (namely, the
separating subspace $S(D)$ of a derivation $D$ of a
(non-)commutative Fr\'{e}chet algebra $A$ is contained in
$\textrm{rad}A$, the (Jacobson) radical of $A$), a stronger
version (namely, the Kaplansky's conjecture $H^1(A, \,A)\,=\,0$;
see also ~\cite{8, 17}), and/or extended versions in
(non-)commutative Fr\'{e}chet algebras, discussed below. For
example, the Singer-Wermer conjecture (and so, its stronger
version, too) does not hold for the derivation
$f\,\rightarrow\,f^{'}$ on the algebra $\textrm{Hol}(U)$, $U$ a
domain in $\cc$, or on the algebra $C^\infty(\rr)$ of smooth
functions on $\rr$ (see also ~\cite{8, 17} for more examples), but
its weaker version does hold since every derivation on a
commutative, semisimple Fr\'{e}chet algebra is continuous by
~\cite{4} and so, $S(D)\,=\,0\,=\,\textrm{rad}A$.

Thus the situation on Fr\'{e}chet algebras is markedly different
from that on Banach algebras, and that a structure theory for
Fr\'{e}chet algebras behaves in a very distinctive manner from
Banach algebras; for example, the prestigious Michael's problem
still remains unsolved since 1952. In fact, this current work may
be considered as a gateway to an affirmative solution to Michael's
problem as follows. The author is finalizing the manuscript in
which the second approach is based on constructing another
Fr\'{e}chet algebra topology (inequivalent to the usual one) on
the test case for this problem, the commutative, semisimple
Fr\'{e}chet algebra $\textrm{Hol}(\ell^\infty)$ of all entire
functions on $\ell^\infty$ ~\cite[\S 3]{13}, using the
discontinuous character (this discontinuous character induces a
discontinuous linear functional on the closed linear subspace of
this algebra which is the first member in the sequence of linear
functionals, giving a sequence of submultiplicative seminorms on
$\textrm{Hol}(\ell^\infty)$; one requires to use the concept of
``tensor product by rows" from ~\cite{13} during the procedure), a
contradiction to the fact that $\textrm{Hol}(\ell^\infty)$ has a
unique Fr\'{e}chet algebra topology as a commutative, semisimple
Fr\'{e}chet algebra of power series (in $\F_\infty$) ~\cite{3, 10,
12}; also, the natural derivation $\partial/\partial X_0$ on the
test case $\textrm{Hol}(\ell^\infty)$ is ``extremely''
discontinuous (i.e., $S(\partial/\partial
X_0)\,=\,\textrm{Hol}(\ell^\infty)$) with respect to another
Fr\'{e}chet algebra topology (cf., Theorem 2.3 below), a
contradiction to the fact that every derivation on the test case
is automatically continuous, by a remark followed by Proposition
2.1 below (cf., also, ~\cite{4}). Some works in the direction of
the non-commutative Singer-Wermer conjecture in the Banach algebra
case have been explored in ~\cite{2, 15}; however, as discussed in
Remark 2.4 of Section 2 below, the non-commutative Singer-Wermer
conjecture does not hold in the Fr\'{e}chet algebra case.

Second, it is easy to give uncountably many inequivalent
Fr\'{e}chet space topologies to a very familiar Fr\'{e}chet space,
namely, spaces of holomorphic functions (see ~\cite{16} for more
details). However the same space is also a semisimple Fr\'{e}chet
algebra of power series, and so, it admits a unique Fr\'{e}chet
algebra topology (see ~\cite{12} for more details). Thus it is
expected that the case of having infinitely many inequivalent
Fr\'{e}chet algebra topologies by some Fr\'{e}chet algebra should
be rare, and the former situation should be common.

We note that in ~\cite{9}, Loy gave a method for constructing
commutative Banach algebras with inequivalent complete norm
topology. We exploit the idea of using the discontinuity of
derivations to give a Fr\'{e}chet algebra other inequivalent
Fr\'{e}chet algebra topologies. We remark that although the
Fr\'{e}chet algebra topology $\tau_0 \;+\;\tau_0$ (defined below)
of $\F_\infty\;\oplus\;\F_\infty$ is not obtainable by our
approach, other inequivalent Fr\'{e}chet algebra topologies on
$\F_\infty\;\oplus\;\F_\infty$ may be constructed. The purpose of
this paper is to show that commutative Fr\'{e}chet algebras
$\F_\infty$ and $\F_\infty\;\oplus\;\F_\infty$ admit countably
many mutually inequivalent Fr\'{e}chet algebra topologies, and to
discuss a couple of avatars of the Singer-Wermer conjecture in
view of these two examples.

Third, we do not know any non-trivial examples of (commutative)
Banach algebras with infinitely many inequivalent complete algebra
norms till this date (see Conclusions below). Of course, one may
consider the trivial product on a Banach space with countably many
inequivalent complete algebra norms to have such an example in the
Banach algebra case.
\section{First Example} In ~\cite{13}, Read showed that the algebra
$\F_\infty \,=\,\cc[[X_0,\,X_1,\, \dots]]$ of all formal power
series in infinitely many commuting indeterminates $X_0,\,X_1,\,
\dots$ has two inequivalent Fr\'{e}chet algebra topologies
$\tau_0$ and $\tau_c$. He also showed that the natural derivation
$\partial/\partial X_0$ is discontinuous with respect to the
Fr\'{e}chet algebra topology $\tau_0$ and with the image the whole
algebra, and that $X_0$ lies in the closure of the coefficient
algebra $\A_0\;=\;\cc[[X_1,\,X_2,\,\dots]]$. Thus, Read's approach
relates constructing ``extremely discontinuous" derivation with
the problem of constructing bizarre Fr\'{e}chet algebra topology
on $\F_\infty$. The motivation for the study of discontinuous
derivations was given by the problem of constructing the dynamics
in statistical mechanics; for more details on the connection
between (discontinuous) derivations and quantum physics, see
Introduction of ~\cite{8}. In this section, we extend Read's
construction. In fact, it is surprising to observe that for each
$i\,\in\,\nn$, the other natural derivation $\partial/\partial
X_i$ is discontinuous with respect to the Fr\'{e}chet algebra
topology $\tau_i$ (defined in the $4^{\textrm{th}}$ paragraph
below), and that for each $i\,\in\,\nn$, $X_i$ lies in the closure
of the coefficient algebra $\A_i\;=\;\cc[[X_0,\,X_1,\,\dots]]$
such that $X_i$ does not belong to $\A_i$. Indeed, with respect to
$\tau_i$, one has $X_n\;\rightarrow \;X_i$ as
$n\;\rightarrow\;\infty$; so one has $X_n - X_i\;\rightarrow \; 0$
yet $\partial/\partial X_i(X_i - X_n)\;=\;1$. Thus $1$ is in the
separating subspace for $\partial/\partial X_i$; since the
subspace is an ideal, it is the whole algebra. This shows that for
each $i\;\in\;\nn$, $\partial/\partial X_i$ vanishes on a dense
subset $\A_i$ of $(\F_\infty,\,\tau_i)$. Not only this, but it is
easy to see that for $i\,\neq\,j$, the Fr\'{e}chet algebra
topologies $\tau_i$ and $\tau_j$ on $\F_\infty$ are not mutually
comparable, because the identity map is continuous in neither
direction as $X_n\,\rightarrow\,X_i$ in $(\F_\infty,\,\tau_i)$
whereas $X_n\,\rightarrow\,X_j$ in $(\F_\infty,\,\tau_j)$. Hence
$\F_\infty$ admits countably many mutually inequivalent
Fr\'{e}chet algebra topologies. Moreover, for $i\,\neq\,j$, the
natural derivations $\partial/\partial X_i$ are continuous on
$(\F_\infty,\,\tau_j)$ as $X_n\,\rightarrow\,X_j$ and
$\partial/\partial X_i(X_j - X_n)\;=\;0$.

In ~\cite{13}, Read claimed that it is not hard to define
continuous derivations on the algebras $\F_\infty$ and $\B$; one
gives them a sensible topology and uses the derivation
$\partial/\partial X_0$ of formal differentiation with respect to
$X_0$. Thus, it was expected that the natural derivation
$\partial/\partial X_0$ is continuous with respect to the
Fr\'{e}chet algebra topology $\tau_c$. By extending Proposition
4.6.3 of ~\cite{5}, we now establish this result in the following
\begin{prop} \label{Prop.4.6.3_Dales} Let $(A,\,\tau)$
be a Fr\'{e}chet algebra containing the indeterminates
$X_0,\,X_1,\,\dots$ and such that $A\,\subset\,\F_\infty$. Suppose
that the embedding $(A, \,\tau)\,\ra\,(\F_n,\,\tau_c)$ is
continuous for some $n\,\in\,\nn$. Then each derivation
$D\,:\,(A,\,\tau)\,\ra\,(A, \,\tau)$ is continuous. In particular,
each derivation $D\,:\,(\F_\infty,\,\tau_c)\,\ra\,(\F_\infty,
\,\tau_c)$ is continuous; more specifically, for each
$i\,\in\,\zz^+$, the natural derivation $\partial/\partial X_i$ on
$(\F_\infty, \,\tau_c)$ is continuous.
\end{prop}
\noindent {\it Proof.} The proof of the first part follows from
proof of Proposition 4.6.3 of ~\cite{5}. For the second part, by
~\cite[Theorem 9.1]{6}, there is a continuous embedding of
$(\F_\infty,\,\tau_c)$ in $(\F_2,\,\tau_c)$. Then, by the first
part, each derivation $D\,:\,\F_\infty\,\ra\,\F_\infty$ is
continuous.
$\hfill \Box$\\

We remark that since the test case for Michael's problem,
discussed above, satisfies the hypotheses of Proposition 2.1 (see,
also, ~\cite[Theorem 10.1]{6}), the above proposition gives
another (elegant) approach to establish the continuity of
derivations on the test case (cf., also, ~\cite{4}). Also, every
derivation on the "natural'' extension of the Beurling-Fr\'{e}chet
algebra of semiweight-type (i.e.,
$\ell^1((\zz^+)^{<\omega},\,\Omega))$ is automatically continuous,
since the Fr\'{e}chet algebra topology $\tau$ on
$\ell^1((\zz+)^{<\omega},\,\Omega)\;(\neq\,\F_\infty)$ is strictly
finer than the Fr\'{e}chet algebra topology $\tau_c$ on
$\F_\infty$ (cf. ~\cite[(4), p. 38]{12}), and so, Theorem 9.1 of
~\cite{6} also holds for this algebra (for the
Beurling-Fr\'{e}chet algebra $\ell^1((\zz^+)^{n},\,\Omega)$ of
semiweight-type, see ~\cite[\S 3, 4]{12}).

We shall feel free to use the terminology and conventions
established in ~\cite{13}. Not only this, but our argument for
this section is kept short because it uses the key ideas involved
in defining the Fr\'{e}chet algebra topology $\tau_0$ on
$\F_\infty$. Let us therefore define ourselves bizarre topology
$\tau_i$ on $\F_\infty$ for each $i\,\in\,\nn$. For this, first,
let us choose a linear functional
$\psi\,:\,A^{(1)}\,\rightarrow\,\cc$ with the property that
$\psi(X_0)\,=\,0$ and $\psi\,=\,X_n^*$ on $A^{(1,
\,0)}\,=\,\B^{(1, \,0)}$ (the subspace of $A^{(1)}$, consisting of
all formal sums $\sum_{j\,=\,1}^{\infty}\lambda_jX_j$), i.e., the
coordinate functional such that
$X_n^*(\sum_{j\,=\,1}^{\infty}\lambda_jX_j)\,=\,\lambda_n$ (and
hence $\psi(X_n)\,=\,1$ for all $n\,\in\,\nn$). Note that this
requires extending the sequence $(X_n)$, $n \in \nn_0$, to a Hamel
basis of $A^{(1)}$, and thus involves the axiom of choice. Then
let $(\phi_n)$, $n \in \nn_0$, be linear functionals on $A^{(1)}$
as follows: (a) For $n \in \nn_0,\;n = i > 0$, $\phi_i = \psi$,
the ``discontinuous'' linear functional, defined above. (b) For
$n\,\in\,\nn_0,\;n\,\neq\,i$, $\phi_n\,=\,X_n^*$, the coordinate
functional such that
$X_n^*(\sum_{j\,=\,0}^{\infty}\lambda_jX_j)\,=\,\lambda_n$. For
each $i
> 0$, let $\tau_i$ be the locally multiplicative convex topology
on $\B$ (a non-commutative analogue of $\F_\infty$, i.e., the
algebra of all formal power series in the infinitely many
non-commutative indeterminates $X_i$) of convergence in all the
seminorms $\|\cdot\|_n^{(m)}$, where
$$\|{\bf a}\|_n^{(m)}\,=\,|{\bf
a}^{(0)}|\,+\,\sum_{r=1}^m\,\sum_{{\bf i}\,\in\,\{0, 1, \dots,
n\}^r}|\otimes_{j=1}^r\phi_{i_j}({\bf a}^{(r)})|$$ for all ${\bf
a}\,\in\,\B$, for this sequence of linear functionals $\phi_n$
(for $i \neq j$, the Fr\'{e}chet algebra topologies $\tau_i$ and
$\tau_j$ are mutually inequivalent, as described above). We remark
that our discontinuous linear functional $\psi$ here is different
from Read's discontinuous linear functional (see ~\cite[Definition
1.7]{13}). In fact, the modification in the definition of $\psi$
is needed so as to change the position of $\phi_i = \psi$ in the
sequence $(\phi_n)$ of linear functionals, in order to ultimately
generate mutually inequivalent Fr\'{e}chet algebra topologies
$\tau_i$ on $\F_\infty$.

We claim that $\B$ is, in fact, complete under the topology
$\tau_i$ for each $i\,>\,0$, so $(\B,\,\tau_i)$ is a Fr\'{e}chet
algebra. Since the averaging map $\alpha\,:\B\,\rightarrow\,\B$,
the locally finite map such that $\alpha(X_{i_1}
\otimes\dots\,\otimes X_{i_n})\,=\,\frac{1}{n!}\sum_{\sigma\in
S_n} X_{i_{\sigma(1)}} \otimes\dots\,\otimes X_{i_{\sigma(n)}}$,
is a $\tau_i$-continuous projection, the subspace $\A = \alpha(\B)
= \ker(I - \alpha)$ is closed, so $(\A,\,\tau_i)$ is also a
Fr\'{e}chet algebra. Note that if we prove our claim we have our
result; for, since $\phi_i(X_n - X_i)\,=\,0$ for $N\,>\,0$, one
sees that $\|X_N - X_i\|_n^{(m)} = 0$ for all $n, m$ with $n < N$.
Hence, $X_N \rightarrow X_i$ in $\tau_i$. One may then use the
derivation $\partial/\partial X_i$ as above, and the separating
subspace is the whole algebra, be it $\F_\infty$ or $\B$. All that
remains is to prove $\B$ complete.

First, we note that $(\B, \,\tau_c)$ is a Fr\'{e}chet algebra; to
define $\tau_c$, apply the method of Corollary 1.11 of ~\cite{13}
to the sequence of coordinate functionals $(X_n^*)$, obtaining
algebra seminorms $|\cdot|_n^{(m)}$, $$|\sum a_{{\bf i}}X^{\otimes
{\bf i}}|_n^{(m)}\,=\,|a^{(0)}|\,+\,\sum_{r=1}^m\,\sum_{{\bf
i}\,\in\,\{0, 1, \dots, n\}^r}|a_{{\bf i}}|.$$ To prove $\B$
complete, by using the seminorms $\|\cdot\|_n^{(m)}$, defined
above, we show that the linear map $\Psi\,:\,\B\,\rightarrow\,\B$,
defined by $$\Psi({\bf a}) = a^{(0)} +
\sum_{r=1}^{\infty}\sum_{{\bf i}\in \nn_0^r}{\bf X}^{\otimes {\bf
i}}\otimes_{j=1}^r\phi_{i_j}({\bf a}^{(r)}),$$ is bijective. (We
note that, for each $i\,\in\,\nn$,
$\Psi\,:\,(\B,\,\tau_i)\,\rightarrow\,(\B,\,\tau_c)$ is
continuous, because convergence under $\tau_i$ is precisely
convergence of all the functionals $\otimes_{j=1}^r\phi_{i_j}({\bf
a}^{(r)})$ that are involved in $\Psi({\bf a})$.) For each $r$,
$\Psi$ maps $\B^{(r)}$ to $\B^{(r)}$, where $\B^{(r)}$ is the
subspace of $n$-homogeneous formal power series $\sum_{{\bf
i}\,\in\,\nn_{0}^{n}}b_{{\bf i}}X^{\otimes\,{\bf i}}$. It is
therefore enough to show that $\Psi^{(r)}\,=\,\Psi|_{\B^{(r)}}$ is
a bijection $\B^{(r)}\,\rightarrow\,\B^{(r)}$ for each $r$. Now,
for $0\,\leq\,k\,\leq\,r$, let $\F^{(r, k)}$ be the set of all
${\bf i}\,=\,(i_1, i_2, \dots, i_r)\,\in\,\nn_0^r$ such that
exactly $k$ of the $i_j$ are equal to zero and let $\B^{(r, k)}$
be the subspace consisting of all power series of the form
$\sum_{{\bf i}\in \F^{(r, k)}}a_{\bf i}X^{\otimes {\bf i}}$. Then
$\B^{(r)}\,=\,\oplus_{k=0}^{r}\B^{(r, \,k)}$. Let $\beta_{r,\,k}$
be the projection onto $\B^{(r, \,k)}$ parallel to the others. We
claim that for each $k$, $\beta_{r,\,k}\Psi^{(r)}|_{\B^{(r,\,k)}}$
is equal to the identity map on $\B^{(r,\,k)}$; and we claim that
$\Psi^{(r)}$ maps $\B^{(r,\,k)}$ into $\oplus_{l=k}^{r}\B^{(r,
\,l)}$. Thus we shall show that the action of $\Psi^{(r)}$ is
``lower triangular'' and nonsingular, with respect to the
decomposition $\B^{(r)}\,=\,\oplus_{k=0}^{r}\B^{(r, \,k)}$.
Therefore, $\Psi^{(r)}$ - and hence $\Psi$ - is bijective.

In view of the definition of $\Psi$ above,
$\beta_{r,\,k}\Psi^{(r)}|_{\B^{(r,\,k)}}$ is equal to the identity
map if and only if we have
$$\otimes_{i=1}^r\phi_{l_i}({\bf a)}\,=\,a_{{\bf
l}}..........(A)$$ for each ${\bf
l}\,=\,(l_1,\,\dots,\,l_r)\,\in\,\F^{(r,\,k)}$ and ${\bf
a}\,\in\,\B^{(r, \,k)}$; and $\Psi^{(r)}$ maps $\B^{(r,\,k)}$ into
$\oplus_{l=k}^{r}\B^{(r, \,l)}$ if and only if
$$\otimes_{i=1}^r\phi_{l_i}({\bf a)}\,=\,0..........(B)$$ for all
$m\,<\,k$, ${\bf l}\,=\,(l_1,\,\dots,\,l_r)\,\in\,\F^{(r,\,m)}$
and ${\bf a}\,\in\,\B^{(r, \,k)}$. Let us prove (A) and (B)
together, by induction on $r$.

When $r\,=\,1$, statement (A) demands that $\phi_{l_1}({\bf
a})\,=\,a_{l_1}$ if either $l_1\,=\,0$ and ${\bf
a}\,\in\,\B^{(1,\,1)}$, or $l_1\,>\,0$ and ${\bf
a}\,\in\,\B^{(1,\,0)}$. For $l_1\,>\,0$, we have two cases: (i) if
$l_1\,\neq\,i$, the functional $\phi_{l_1}\,=\,X_{l_1}^*$ anyway,
so the assertion is true; and (ii) if $l_1\,=\,i$,
$\phi_i\,=\,\psi\,=\,X_i^*$ on $\B^{(1,\,0)}$, so the assertion is
true. But when $l_1\,=\,0$, we note that $\B^{(1,\,1)}$ consists
solely of multiples $a_0X_0$; since $\phi_0\,=\,X_0^*$, the
assertion is true for $l_1\,=\,0$ also. Statement (B) demands that
$\phi_{l_1}({\bf a})\,=\,0$ if $l_1\,>\,0$ and ${\bf
a}\,\in\,\B^{(1,\,1)}$. Again we have two cases: (i) if
$l_1\,\neq\,i$, the functional $\phi_{l_1}\,=\,X_{l_1}^*$, so
$X_{l_1}^*(a_0X_0)\,=\,0$; and (ii) if $l_1\,=\,i$,
$\phi_i\,=\,\psi$, so $\psi(a_0X_0)\,=\,0$. So the assertion is
correct.

When $r\,>\,1$, we proceed by induction. For this, we need the
following definition, taken from ~\cite{13}.\\

\begin{defn} \label{Definition 2.2.} For each $n\,\in\,\nn_0$, let
$P_n\,:\,\B\,\rightarrow\,\B$ be the locally finite map such that
$P_n(1)\,=\,0$ and
\begin{equation}
P_n(X_{i_{1}}\,\otimes\,X_{i_{2}}\,\otimes\,\cdots\,\otimes\,X_{i_{m}})
= \begin{cases} 0, ~&\text{when}\;i_{1}\,\neq\,n,\\
X_{i_{2}}\,\otimes\,\cdots\,\otimes\,X_{i_{m}}, ~&\text{when}\;
i_{1}\,=\,n. \end{cases}
\end{equation}
Thus $P_n$ takes the quotient on division from the left by $X_n$,
and discards the remainder.
\end{defn}

First we establish statements (A) and (B) when $r\,=\,2$. In this
case, from Definition 1.9 of ~\cite{13}, we have the ``tensor
product by rows" $\phi_{l_1}\,\otimes
\,\phi_{l_2}\,:\,\B^{(2)}\,\rightarrow \,\cc$ by
$$\phi_{l_1}\,\otimes \,\phi_{l_2}({\bf
a})\,=\,\phi_{l_1}(\sum_{j=0}^{\infty}X_j\phi_{l_2}(P_j{\bf
a))}..........(C)$$ To prove statement (A), consider the case when
${\bf l}\,=\,(l_1,\,l_2)\,\in\,\F^{(2,\,k)}$ and ${\bf
a}\,\in\,\B^{(2, \,k)}$. Note that $P_j{\bf a}$ lies in $\B^{(1,
\,1)}$ if $j\,>\,0$, but in $\B^{(1, \,0)}$ if $j\,=\,0$, because
the division on the left by $X_0$ removes a factor $X_0$ from each
monomial ${\bf X^{\otimes i}}$. If $l_1\,=\,0$, then
$(l_2)\,\in\,\F^{(1,\,0)}$ (i.e., $l_2\,>\,0$), so
$\phi_{l_2}({\bf a})\,=\,0$ when ${\bf a}\,\in\,\B^{(1, \,1)}$
(here, if $l_2\,=\,i$, then $\phi_i\,=\,\psi$, so
$\psi(a_0X_0)\,=\,0$; otherwise, $\phi_{l_2}\,=\,X_{l_2}^*$ and
$X_{l_2}^*(a_0X_0)\,=\,0$), but will send ${\bf
b}\,\in\,\B^{(1,\,0)}$ to $b_{l_2}$ (here, ${\bf
b}\,=\,\sum_{j=1}^{\infty}b_jX_j$ and if $l_2\,=\,i$, then
$\phi_i\,=\,\psi\,=\,X_i^*$ on $\B^{(1,\,0)}$; otherwise,
$\phi_{l_2}\,=\,X_{l_2}^{*}$ and
$\phi_{l_2}(\sum_{j\,=\,1}^{\infty}b_jX_j)\,=\,b_{l_2}$).
Accordingly, when $l_1\,=\,0$, statement (C) is equal to
$\phi_0(X_0(P_0{\bf a})_{l_2})\,=\,(P_0{\bf
a})_{l_2}\,=\,a_{l_1,\,l_2}$. If on the other hand, $l_1\,>\,0$,
then we have two cases: (i) if $l_1\,=\,i$, then
$\phi_{l_1}\,=\,\phi_i$. So, the statement (C) is equal to
$\phi_i(\sum_{j=0}^{\infty}X_j\phi_{l_2}(P_j{\bf
a}))\,=\,\phi_i(X_0\phi_{l_2}(P_0{\bf a}))
\,+\,\phi_i(\sum_{j=1}^{\infty}X_j\phi_{l_2}(P_j{\bf
a}))\,=\,\phi_i(\sum_{j=1}^{\infty}X_j\phi_{l_2}(P_j{\bf
a}))\,=\,\phi_{l_2}(P_i{\bf a}))\,=\,\phi_0(a_0X_0)$ as
$\phi_i\,=\,\psi\,=\,X_i^*$ on $\B^{(1,\,0)}$. In this case,
$l_2\,\in\,\F^{(1,\,1)}$ and $P_i{\bf a}\,\in\,\B^{(1,\,1)}$. So,
the statement (C) is equal to $(P_i{\bf
a})_{l_2}\,=\,a_{l_1,\,l_2}$ as before. (ii) if $l_1\,\neq\,i$,
then $\phi_{l_1}\,=\,X_{l_1}^*$, so the statement (C) is equal to
$\phi_{l_2}(P_{l_1}{\bf a})$. In this case, $(l_2) \,\in\,\F^{(1,
k)}$ and $P_{l_1}{\bf a}\,\in\,\B^{(1, k)}$ so the statement (C)
is equal to $(P_{l_1}{\bf a})_{l_2}\,=\,a_{l_1,\, l_2}$ as before.

For the statement (B), consider ${\bf a}\,\in\,\B^{(2, \,k)}$ as
before, but now use sequences ${\bf l}\,\in\,\F^{(r,\,m)}$ for
some $m\,<\,k$. If $l_1\,=\,0$, then $l_2\,\in\,\F^{(1,\,0)}$;
since each $P_j{\bf a}$ is either in $\B^{(1,\,1)}$ or
$\B^{(1,\,0)}$, $\phi_{l_2}(P_j{\bf a})\,=\,0$ for all $P_j{\bf
a}$ (even when $l_2\,=\,i$). If $l_1\,>\,0$, then we have two
cases: (i) if $l_1\,=\,i$, then $\phi_{l_1}\,=\,\phi_i\,=\,\psi$,
so $\phi_i(\sum_{j=0}^{\infty}X_j\phi_{l_2}(P_j{\bf
a}))\,=\,\phi_{l_2}(P_i{\bf a})$. Since $l_2\,\in\,\F^{(1,\,m)}$
and $P_{l_1}{\bf a}\,\in\,\B^{(1,\,1)}$, $\phi_{l_2}(P_{l_1}{\bf
a})\,=\,0$. (ii) if $l_1\,\neq\,i$, then
$\phi_{l_1}\,=\,X_{l_1}^*$, so the statement (C) is equal to
$\phi_{l_2}(P_{l_1}{\bf a})$. Since $(l_2) \,\in\,\F^{(1, m)}$ and
$P_{l_1}{\bf a}\,\in\,\B^{(1, k)}$ so this is zero also. Thus the
statement (A) and (B) are proved for $r = 2$.

Now one establishes statements (A) and (B) by induction when
$r\,>\,2$. For this, follow Read's argument in our case.
\begin{thm} \label{Theroem 2.5_Read} Let $i \in \zz^+$ be fixed.
$(\B,\,\tau_i)$ is complete with respect to the seminorms
$\|\cdot\|_n^{(m)}$. The derivation $\partial/\partial X_i :
(\B,\,\tau_i) \rightarrow (\B,\,\tau_i)$ is discontinuous, and its
separating subspace is all of $\B$. The derivation
$\partial/\partial X_i : (\A,\,\tau_i) \rightarrow (\A,\,\tau_i)$
is also discontinuous, and its separating subspace is all of $\A$.
\end{thm}
\noindent {\it Proof.} Let $(a_n)$ be a Cauchy sequence in
$(\B,\,\tau_i)$. Then $\Psi(a_n)$ is a Cauchy sequence in
$(\B,\,\tau_c)$. Since $\B$ is complete with respect to the
topology $\tau_c$, let $a$ be the $\tau_c$-limit of $(\Psi(a_n))$.
The $\tau_i$-limit of $(a_n)$ is just $\Psi^{-1}(a)$. So
$(\B,\,\tau_i)$ is complete.

The argument that shows the rest of the assertions, given $(\B,
\,\tau_i)$ complete, has already been given above. $\hfill \Box$\\

\noindent {\bf Remark 2.4.} The above theorem shows that the
non-commutative analogues of all avatars of the Singer-Wermer
conjecture do not hold in the Fr\'{e}chet case.
\section{Second Example} We now show that the algebra $\F_\infty\;\oplus\;\F_\infty$ has
countably many mutually inequivalent Fr\'{e}chet algebra
topologies. In fact, the algebra $\F_\infty$ has been completed by
the adjunction of a radical so that $\F_\infty\;\oplus\;\F_\infty$
has Fr\'{e}chet algebra topologies in which $\F_\infty$ is dense
in those topologies.

For us, let $A = \F_\infty$ be a commutative (Fr\'{e}chet) algebra
and $M \;=\;\F_\infty$ a commutative (Fr\'{e}chet) $A$-module. Let
$H^1(\F_\infty,\;\F_\infty)$ denote the first algebraic cohomology
group, $H_C^1(\F_\infty,\;\F_\infty)$ the first continuous
cohomology group, where the cochains are required to be bounded.
Thus with the usual conventions $H^1(\F_\infty,\;\F_\infty)$ is
the space of derivations of $\F_\infty$ into itself;
$H_C^1(\F_\infty,\;\F_\infty)$ the space of continuous derivations
of $\F_\infty$ into itself with the Fr\'{e}chet algebra topology
$\tau_c$ on $\F_\infty$.

Let $\F_\infty$ be a commutative Fr\'{e}chet algebra with
Fr\'{e}chet algebra topology $\tau_i$ generated by the sequence
$(p_{k, \,i}^{'})$ for each $i\,\in\,\zz^+$, by \S 2 and $\tau_c$
generated by the sequence $(p_k)$. Let $\A$ denote the vector
space direct sum $\F_\infty\;\oplus\;\F_\infty$ with product
$$(x, m) (y, n) \,=\,(xy, x\cdot n + y\cdot m),$$ seminorms
$$q_k(x, m)\,=\,p_k(x) + p_k(m),$$ generating the Fr\'{e}chet algebra topology
$\tau_c\;+\;\tau_c$, and seminorms $$q_{k, \,i}^{'}(x,
m)\,=\,p_{k, \,i}^{'}(x) + p_{k, \,i}^{'}(m),$$ generating the
Fr\'{e}chet algebra topology $\tau_i\;+\;\tau_i$ for each
$i\,\in\,\zz^+$. For each $i\;\in\;\zz^+$,
$\partial_i\,=\,\partial/\partial
X_i\,\in\,H^1(\F_\infty,\;\F_\infty)$, the functional
$$q_{k, \partial_i}\,:\,(x, m)\,\rightarrow\,p_k(x) + p_k(\partial_i(x)-m)$$ is
defined on the algebra $\F_\infty\;\oplus\;\F_\infty$ and is, in
fact, easily seen to be a submultiplicative seminorm thereon; for
each $i\;\in\;\zz^+$, the sequence $(q_{k, \partial_i})$ generates
the Fr\'{e}chet algebra topology $\tau_{\partial_i}$ equivalent to
the Fr\'{e}chet algebra topology $\tau_c\;+\;\tau_c$ (see Theorem
~\ref{Thm. 5.1.17_Frechet} below). For each $i\;\in\;\zz^+$,
$\partial_i\,=\,\partial/\partial
X_i\,\in\,H^1(\F_\infty,\;\F_\infty)$, the functional
$$q_{k, \partial_i}^{'}\,:\,(x, m)\,\rightarrow\,p_{k, \,i}^{'}(x) + p_{k, \,i}^{'}(\partial_i(x)-m)$$ is
defined on the algebra $\F_\infty \oplus \F_\infty$ and is, in
fact, easily seen to be a submultiplicative seminorm thereon; for
each $i\;\in\;\zz^+$, the sequence $(q_{k, \partial_i}^{'})$
generates the Fr\'{e}chet algebra topology $\tau_{\partial_i}^{'}$
inequivalent to the Fr\'{e}chet algebra topology
$\tau_i\;+\;\tau_i$ (see Theorem ~\ref{Thm. 5.1.17_Frechet}
below). For each $i\;\in\;\zz^+$,
$\partial_i\,=\,\partial/\partial X_i$, the map
$$\theta_{\partial_{i}}\,:\,(x, m)\,\rightarrow\,(x,\; \partial_i(x) - m)$$ is an isometric
isomorphism of $\F_\infty \oplus \F_\infty$ under $(q_{k,
\partial_i})$ (respectively, $(q_{k, \partial_i}^{'})$) into
$(\F_\infty \oplus \F_\infty,\;\tau_c + \tau_c)$ (respectively,
$(\F_\infty \oplus \F_\infty,\;\tau_i + \tau_i)$) and so, extends
uniquely to a map of the completion $(\F_\infty \oplus
\F_\infty)_{\partial_i}$ of $\F_\infty \oplus \F_\infty$ into
$(\F_\infty \oplus \F_\infty,\;\tau_c + \tau_c)$ (respectively,
$(\F_\infty \oplus \F_\infty,\;\tau_i + \tau_i)$). In particular,
if $\iota\,:\,x\,\rightarrow\,(x, 0)$ is the natural embedding of
$\F_\infty$ into $\F_\infty \oplus \F_\infty$, then, for each
$i\;\in\;\zz^+$, $q_{k,
\partial_i}^{''}\,:\,x\,\rightarrow\,q_{k, \partial_i}(\iota(x))$ (respectively,
$q_{k, \partial_i}^{'''}\,:\,x\,\rightarrow\,q_{k,
\partial_i}^{'}(\iota(x))$) is a seminorm on $\F_\infty$ and $\theta_{\partial_i}
\,\circ\,\iota$ extends to an isometric isomorphism of
$\overline{(\F_\infty)}_{\partial_i}$, the completion of
$\F_\infty$ under $(q_{k, \partial_i}^{''})$ (respectively,
$(q_{k, \partial_i}^{'''})$), with
$\overline{\textrm{Gr}(\partial_i)}$, the closure (respectively,
in $(\F_\infty \oplus \F_\infty,\;\tau_c + \tau_c)$ and in
$(\F_\infty \oplus \F_\infty,\;\tau_i + \tau_i)$) of the graph of
$\partial_i$.

Let $i\,\in\,\zz^+$ be fixed. Now if $\partial_i$ is continuous on
$(\F_\infty, \,\tau_c)$, then $(q_{k,
\partial_i})$ is equivalent to $(q_k)$ on
$\F_\infty \oplus \F_\infty$ (see Theorem ~\ref{Thm.
5.1.17_Frechet} below) and $(q_{k,
\partial_i}^{''})$ is equivalent to $(p_k)$ on $\F_\infty$. Thus
$(\F_\infty \oplus \F_\infty)_{\partial_i}\,=\,(\F_\infty \oplus
\F_\infty,\;\tau_c + \tau_c)$ and
$\overline{(\F_\infty)}_{\partial_i}\,=\,(\F_\infty,\;\tau_c)$.
This is the case when $\F_\infty$ is a Fr\'{e}chet algebra under
the Fr\'{e}chet algebra topology $\tau_c$. In the discontinuous
case $q_{k, \partial_i}^{'''}$ is a discontinuous seminorm on
$(\F_\infty,\;\tau_i)$ and $\iota$ is a discontinuous isomorphism.
We study this case in detail. First, by Theorem ~\ref{Thm.
5.1.17_Frechet} below, $(q_{k, \partial_i}^{'})$ is not equivalent
to $(q_{k, \,i}^{'})$ on $\F_\infty\;\oplus\;\F_\infty$.

Suppose now that $\partial_i$ is discontinuous on $(\F_\infty,
\,\tau_i)$. Then we have
$$\overline{\textrm{Gr}(\partial_i)}^{\tau_i\;+\tau_i}\,\bigcap\,0\,\oplus\,\F_\infty\,\neq\,\{0\}.$$
Thus, if $\F_\infty$ is a Fr\'{e}chet algebra with
$H^1(\F_\infty,\;\F_\infty)\,\neq\,H_C^1(\F_\infty,\;\F_\infty)$,
then $\F_\infty$ has a completion with a non-trivial nil ideal,
that is, the completion $\overline{(\F_\infty)}_{\partial_i}$
under $(q_{k, \partial_i}^{'''})$ is, in fact,
$(\F_\infty\;\oplus\;\F_\infty,\;\tau_i\;+\;\tau_i)$ due to the
following theorem. We remark that if $A$ is a commutative,
semisimple Fr\'{e}chet algebra, then $H^1(A, A)\,=\,H_C^1(A, A)$
(see ~\cite{4} for details), and in particular, if $A$ is a
commutative pro-$C^*$-algebra, then $H^1(A, A)\;=\;0$ (see
~\cite{17} for more general results). However we do not know an
example of a non-semisimple, non-Banach Fr\'{e}chet algebra such
that the extended Kaplansky's conjecture that $H^1(A, M)\,=\,0$
for any Fr\'{e}chet $A$-module $M$ holds (in particular, such that
the Kaplansky's conjecture that $H^1(A, A)\,=\,0$ holds). The
Singer-Wermer conjecture holds for commutative Banach algebras;
however, we do not know whether its extended version that the
image of a derivation $D\,:\,A\,\ra\,M$ is contained in the
radical of $M$, where $A$ is a commutative Banach algebra and $M$
is a commutative Banach $A$-module, holds (in the non-normed
settings, we must consider the weaker form of this, i.e.,
$S(D)\,\subset\,\textrm{rad}M$; see Introduction above). Banach
algebras with $H^1(A, M)\,=\,H_C^1(A, M)$ are discussed in
~\cite{1}.
\begin{thm} \label{Theorem 1_Loy} Let $i\,\in\,\zz^+$ be fixed.
Let $(\F_\infty,\;\tau_i)$ be a commutative Fr\'{e}chet algebra,
$\partial_i\;=\;\partial/\partial X_i$ a non-zero derivation of
$\F_\infty$ into itself. Then the algebra
$\overline{(\F_\infty)}_{\partial_i}$ admits a Fr\'{e}chet algebra
topology $\tau_{\partial_i}^{'}$, generated by
$(q_{k,\;\partial_i}^{'})$, and which is inequivalent to the
Fr\'{e}chet algebra topology $\tau_i\;+\;\tau_i$, generated by
$(q_{k, \,i}^{'})$. In particular, the algebra
$\F_\infty\;\oplus\;\F_\infty$ admits countably many Fr\'{e}chet
algebra topologies inequivalent to the Fr\'{e}chet algebra
topology $\tau_i\;+\;\tau_i$.
\end{thm}
\noindent {\it Proof.} First, by \S 2, we note that, for each
$i\;\in\;\zz^+$, $\partial_i$ vanishes on a dense subset $\A_i$
(which is the coefficient algebra of $\F_\infty$); we also note
that this is possible due to the unique property of the
Fr\'{e}chet algebra topology $\tau_i$ as explained in \S 2. Now
follow the proof of ~\cite[Theorem 1]{9} for this case. $\hfill
\Box$\\

\begin{rem} \label{Remark 3.2.} We note that $(\F_\infty,\;\tau_i)$ is
dense in $$\overline{(\F_\infty)}_{\partial_i}=(\F_\infty\oplus
\F_\infty,\;\tau_i + \tau_i),$$ and
$\textrm{Rad}(\F_\infty\;\oplus\;\F_\infty)\;=\;\textrm{Rad}(\F_\infty)\;\oplus\;\F_\infty$.
\end{rem}
\begin{rem} \label{Remark 3.3.} Theorem ~\ref{Theorem 1_Loy} gives
countably many mutually inequivalent Fr\'{e}chet algebra
topologies on $\F_\infty \oplus \F_\infty$; for mutual
inequivalence, see Theorem 3.9 below. In fact, we have a stronger
result than Theorem ~\ref{Theorem 1_Loy} as follows. We note that,
for each $i\;\in\;\zz^+$, the derivation $\partial_i$ on
$\F_\infty$ induces a natural derivation $D_i \;:\;(a,
\;x)\;\mapsto\;(0,\;\partial_i(a))$ on
$\F_\infty\;\oplus\;\F_\infty$.
\end{rem}
\begin{thm} \label{Thm. 5.1.17_Frechet} The algebra $\F_\infty\;\oplus\;\F_\infty$
is a Fr\'{e}chet algebra with respect to the sequences $(q_k)$,
$(q_{k, \,i}^{'})$, $(q_{k, \partial_i})$ and $(q_{k,
\partial_i}^{'})$. The Fr\'{e}chet algebra topologies
$\tau_{\partial_i}$ generated by the sequences
$(q_{k,\partial_i})$ are equivalent to the Fr\'{e}chet algebra
topology $\tau_c\;+\;\tau_c$ if and only if for each
$i\;\in\;\zz^+$, $\partial_i$ is continuous on
$(\F_\infty,\;\tau_c)$ if and only if for each $i\;\in\;\zz^+$,
the natural derivation $D_i$ is continuous on
$(\F_\infty\;\oplus\;\F_\infty,\;\tau_c\;+\;\tau_c)$. Let
$i\,\in\,\zz^+$ be fixed. Then the Fr\'{e}chet algebra topology
$\tau_{\partial_i}^{'}$ generated by the sequence
$(q_{k,\partial_i}^{'})$ is not equivalent to the Fr\'{e}chet
algebra topology $\tau_i\;+\;\tau_i$ if and only if $\partial_i$
is discontinuous on $(\F_\infty,\;\tau_i)$ if and only if the
natural derivation $D_i$ is discontinuous on
$(\F_\infty\;\oplus\;\F_\infty,\;\tau_i\;+\;\tau_i)$.
\end{thm}
\noindent {\it Proof.} Certainly $(\F_\infty\;\oplus\;\F_\infty,
(q_k))$ (respectively, for each $i\;\in\;\zz^+$,
$(\F_\infty\;\oplus\;\F_\infty,\, (q_{k,\,i}^{'}))$) is a
Fr\'{e}chet algebra and $q_{k, \partial_i}$ (respectively, $q_{k,
\partial_i}^{'}$) is a seminorm on $\F_\infty\;\oplus\;\F_\infty$ for each
$k\,\in\,\nn$. For $(a, x),\, (b,
y)\,\in\,\F_\infty\;\oplus\;\F_\infty$, we have
\begin{align*}
q_{k,
\partial_i}((a, x)(b, y)) &= p_k(ab)+p_k(a \cdot (\partial_i(b) -
y) + b
\cdot (\partial_i(a) - x))\\
&\leq (p_k(a) + p_k(\partial_i(a) - x)) (p_k(b) +
p_k(\partial_i(b) - y))\\
&= q_{k, \partial_i}((a, x)) q_{k, \partial_i}((b, y)),
\end{align*} and so $ q_{k,\partial_i}$ (respectively, $q_{k,
\partial_i}^{'}$) is a submultiplicative seminorm on
$\F_\infty\;\oplus\;\F_\infty$ for each $k\,\in\,\nn$. We now show
that $(\F_\infty\;\oplus\;\F_\infty,\, (q_{k,
\partial_i}))$ (respectively, for each $i\;\in\;\zz^+$, $(\F_\infty\;\oplus\;\F_\infty,
(q_{k, \partial_i}^{'}))$) is a Fr\'{e}chet algebra. Let $((a_n,
x_n))$ be a Cauchy sequence in $(\F_\infty\;\oplus\;\F_\infty,\,
(q_{k, \partial_i}))$ (respectively, for each $i\;\in\;\zz^+$,
$(\F_\infty\;\oplus\;\F_\infty, (q_{k, \partial_i}^{'}))$). Then
$(a_n)$ and $(\partial_i(a_n) - x_n)$ are Cauchy sequences in
$(\F_\infty, \,(p_k))$ (respectively, for each $i\;\in\;\zz^+$,
$(\F_\infty, \,(p_{k,\,i}^{'}))$. Since $\F_\infty$ is a
Fr\'{e}chet space, there exists $a\,\in\,\F_\infty$ and
$x\,\in\,\F_\infty$ such that $a_n\,\rightarrow\,a$ and
$\partial_i(a_n) - x_n\,\rightarrow\,x$. Then $(a_n,
x_n)\,\rightarrow\,(a,
\partial_i(a) - x))$ in $(\F_\infty\;\oplus\;\F_\infty,\, (q_{k,
\partial_i}))$ (respectively, for each $i \in \zz^+$, $(\F_\infty \oplus \F_\infty,
(q_{k,
\partial_i}^{'}))$) and so $(\F_\infty\;\oplus\;\F_\infty,\,
(q_{k, \partial_i}))$ (respectively, for each $i\;\in\;\zz^+$,
$(\F_\infty\;\oplus\;\F_\infty, (q_{k, \partial_i}^{'}))$) is a
Fr\'{e}chet algebra.

Suppose that $\partial_i$ is continuous. Then, for each
$m\,\in\,\nn$, there exists $n(m) \in \nn$ and a constant $c_m >
0$ such that $$q_{m,
\partial_i}((a, x))\leq
p_{m}(a)+c_mp_{n(m)}(a)+p_m(x)\leq(1+c_m)q_{n(m)}((a, x)),$$ and
so the two topologies are equivalent, by the open mapping theorem
for Fr\'{e}chet spaces.

Conversely, suppose that the two topologies are equivalent on the
algebra $\F_\infty\;\oplus\;\F_\infty$. Then, for each
$m\,\in\,\nn$, there exists $n(m)\,\in\,\nn$ and a constant $c_m
> 0$ such that $$q_{m, \partial_i}((a,
x)) \leq c_m\,q_{n(m)}((a, x))\;((a, x)\,\in\,\F_\infty \oplus
\F_\infty).$$ Hence, $$p_m(\partial_i(a))\,\leq\, q_{m,
\partial_i}((a, 0))\,\leq\,c_m\,q_{n(m)}((a,
0))\,=\,c_m\,p_{n(m)}(a)\;\;(a\,\in\,\F_\infty),$$ and so,
$\partial_i$ is continuous on $(\F_\infty,\;\tau_c)$.

The first half of the proof of the second statement has already
been discussed in the proof of Theorem 3.1. Suppose that for each
$i \in \zz^+$, the Fr\'{e}chet algebra topology
$\tau_{\partial_i}^{'}$ is not equivalent to the Fr\'{e}chet
algebra topology $\tau_i + \tau_i$ on $\F_\infty \oplus
\F_\infty$. Then, for all $n \in \nn$ and a constant $c > 0$,
there exists $m \in \nn$ such that
$$q_{m,\partial_i}^{'}((a,\;x))\;>\;c\;q_{n,\,i}^{'}((a,\;x))\;((a,\;x)\;\in\;\F_\infty\;\oplus\;\F_\infty).$$
Hence,
$$p_{m,\,i}^{'}(\partial_i(a))\;\geq\;q_{m,\partial_i}^{'}((0,\,a))\;>\;cq_{n,\,i}^{'}((0,
\,a))\;=\;cp_{n,\,i}^{'}(a)\;(a\;\in\;\F_\infty),$$ and so,
$\partial_i$ is discontinuous on $(\F_\infty,\;\tau_i)$. $\hfill
\Box$\\

\begin{rem} \label{Remark 3.5.} For each $i\;\in\;\zz^+$, the natural
derivation $D_i$ does not vanish on a dense subset $\F_\infty$ of
the algebra $(\F_\infty\;\oplus\;\F_\infty,\;\tau_i\;+\;\tau_i)$;
however, it does vanish on the coefficient algebra $\A_i$ since
the restriction of the natural derivation $D_i$ on $\F_\infty$ is
$\partial_i$.
\end{rem}

\begin{rem} \label{Remark 3.6.} It is easy to see that for each
$i\,\in\,\zz^+$, the Fr\'{e}chet algebra topologies
$\tau_c\;+\;\tau_c$ and $\tau_i\;+\;\tau_i$ are inequivalent,
since if they are equivalent on the algebra
$\F_\infty\;\oplus\;\F_\infty$, then since $\F_\infty$ can be
identified with the closed subalgebra $\F_\infty\;\oplus\;0$, we
see that for each $i\,\in\,\zz^+$, the Fr\'{e}chet algebra
topologies $\tau_c$ and $\tau_i$ are equivalent, a contradiction.
For the same reason, for $i\,\in\,\zz^+$, the Fr\'{e}chet algebra
topologies $\tau_i\;+\;\tau_i$ are mutually inequivalent on the
algebra $\F_\infty\;\oplus\;\F_\infty$. Similarly, for
$i\,\in\,\zz^+$, the Fr\'{e}chet algebra topologies
$\tau_c\;+\;\tau_c$ and $\tau_{\partial_i}^{'}$ are inequivalent,
since if they are equivalent on the algebra
$\F_\infty\;\oplus\;\F_\infty$, then, by Theorem ~\ref{Thm.
5.1.17_Frechet} and by the transitivity property of the relation
``equivalence of topologies'', the Fr\'{e}chet algebra topologies
$\tau_{\partial_i}$ and $\tau_{\partial_i}^{'}$ are equivalent on
the algebra $\F_\infty\;\oplus\;\F_\infty$, a contradiction to the
remark, given after Theorem 3.9, below. Also, for each
$i\,\in\,\zz^+$, the Fr\'{e}chet algebra topologies
$\tau_i\;+\;\tau_i$ and $\tau_{\partial_i}$ are inequivalent,
since if they are equivalent on the algebra
$\F_\infty\;\oplus\;\F_\infty$, then, by Theorem ~\ref{Thm.
5.1.17_Frechet} and by the transitivity property of the relation
``equivalence of topologies'', the Fr\'{e}chet algebra topologies
$\tau_i\;+\;\tau_i$ and $\tau_c\;+\;\tau_c$ are equivalent, a
contradiction to the fact given above.
\end{rem}
From the first half of Theorem ~\ref{Thm. 5.1.17_Frechet}, we have
the following
\begin{cor} \label{Thm. 5.1.17_corollary1} The
Fr\'{e}chet algebra topologies $\tau_{\partial_i}$ are mutually
equivalent on the algebra $\F_\infty\;\oplus\;\F_\infty$.
\end{cor}
\noindent {\it Proof.} Clearly, by Theorem ~\ref{Thm.
5.1.17_Frechet}, for each $i\;\in\;\zz^+$, the Fr\'{e}chet algebra
topology $\tau_{\partial_i}$ generated by the sequence
$(q_{k,\partial_i})$ is equivalent to the Fr\'{e}chet algebra
topology $\tau_c\;+\;\tau_c$, since the derivations $\partial_i$
is continuous on $(\F_\infty,\,\tau_c)$. Now, the Fr\'{e}chet
algebra topology $\tau_{\partial_i}$ is equivalent to
$\tau_{\partial_j}$ by the transitivity property of the relation
``equivalence of
topologies''. $\hfill \Box$\\

In fact, we can sharpen the above corollary as follows.
\begin{cor} \label{Thm. 5.1.17_corollary1extended} Let $i$ and $j$ be in $\zz^+$ such that $i\,\neq\,j$.
The Fr\'{e}chet algebra topology $\tau_{\partial_i}$ is equivalent
to the Fr\'{e}chet algebra topology $\tau_{\partial_j}$ on the
algebra $\F_\infty\;\oplus\;\F_\infty$ if and only if the natural
derivation $D_i$ is continuous on
$(\F_\infty\;\oplus\;\F_\infty,\;\tau_{\partial_j})$ and the
natural derivation $D_j$ is continuous on
$(\F_\infty\;\oplus\;\F_\infty,\;\tau_{\partial_i})$.
\end{cor}
\noindent {\it Proof.} Suppose that for $i\;\neq\;j$, the natural
derivation $D_i$ is continuous on
$(\F_\infty\;\oplus\;\F_\infty,\;\tau_{\partial_j})$ and the
natural derivation $D_j$ is continuous on
$(\F_\infty\;\oplus\;\F_\infty,\;\tau_{\partial_i})$. Then, for
each $m\;\in\;\nn$ and $(a,
x)\;\in\;\F_\infty\;\oplus\;\F_\infty$, there exist
$n(m)\;\in\;\nn$ and constants $c_m,\,c_m^{'}\;>\;0$ such that
\begin{align*} q_{m,\partial_i}((a,\;x))&\leq p_m(a) + p_m(\partial_i(a))
+ p_m(x)\\
&\leq p_m(a) + c_mp_{n(m)}(a) + p_m(x)\\ &\leq (1 +
c_m)q_{n(m)}((a, x))\\ &\leq (1 +
c_m)c_m^{'}q_{n(m),\partial_j}((a, x)),
\end{align*}
and so, the two topologies $\tau_{\partial_i}$ and
$\tau_{\partial_j}$ are equivalent on the algebra
$\F_\infty\;\oplus\;\F_\infty$, by the open mapping theorem for
Fr\'{e}chet spaces.

Conversely, suppose that for $i \neq j$, the two topologies
$\tau_{\partial_i}$ and $\tau_{\partial_j}$ are equivalent on the
algebra $\F_\infty \oplus \F_\infty$. Then, for each $m \in \nn$
and $(a, x) \in \F_\infty \oplus \F_\infty$, there exist $n(m),
l(m)\;\in\;\nn$ and constants $c_m, c_m^{'}\;>\;0$ such that
$q_{m,\partial_i}((a, x))\;\leq\;c_mq_{n(m),\partial_j}((a, x)).$
Hence, $$q_{m,\partial_i}(D_j(a,
x))\;\leq\;c_mq_{n(m),\partial_j}(D_j(a,
x))\;\leq\;c_m^{'}q_{l(m)}(D_j(a, x))$$ since the sequences
$(q_{k, \partial_i})$ and $(q_k)$ are equivalent and
$$c_m^{'}q_{l(m)}(D_j(a, x))\leq\;c_m^{'}q_{l(m)}(a, x)\;\leq\;c_m^{'}q_{l(m),\partial_i}(a,
x)$$ since $D_j$ is continuous on
$(\F_\infty\;\oplus\;\F_\infty,\;\tau_c + \tau_c)$. Thus, for
$i\;\neq\;j$, the natural derivation $D_i$ is continuous on
$(\F_\infty\;\oplus\;\F_\infty,\;\tau_{\partial_j})$ and the
natural derivation $D_j$ is continuous on
$(\F_\infty\;\oplus\;\F_\infty,\;\tau_{\partial_i})$. $\hfill
\Box$
\begin{thm} \label{Thm. 5.1.17_corollary} The
Fr\'{e}chet algebra topologies $\tau_{\partial_i}^{'}$ are
mutually inequivalent on the algebra
$\F_\infty\,\oplus\,\F_\infty$ if and only if for each $i,\;
j\;\in\;\zz^+$ such that $i\;\neq\;j$, the natural derivation
$D_i$ is discontinuous on
$(\F_\infty\;\oplus\;\F_\infty,\;\tau_{\partial_j}^{'})$.
\end{thm}
\noindent {\it Proof.} Let
$(\F_\infty\;\oplus\;\F_\infty,\;\tau_{\partial_j}^{'})$ be a
commutative Fr\'{e}chet algebra and let $D_i$ be the discontinuous
natural derivation on $\F_\infty\;\oplus\;\F_\infty$, induced by
$\partial_i$. Then, following the arguments of Theorem 1 of
~\cite{9},
$$S\;=\;\{m\;=\;(m_1, m_2)\;\in\;\F_\infty\;\oplus\;\F_\infty\;:\;(0,m)\;\in\;\overline{\textrm{Gr}(D_i)}^{\tau_{\partial_j}^{'}\;+\tau_{\partial_j}^{'}}\},$$
so that $S$ is certainly contained in
$\overline{D_i(\F_\infty\;\oplus\;\F_\infty)}^{\tau_{\partial_j}^{'}}$.
Now a diagonal argument shows that $S$ is closed. Then, following
the arguments of Theorem 1 of ~\cite{9}, we have
$\F_\infty\;\oplus\;\F_\infty\;\oplus\;S\,=\,\overline{\textrm{Gr}(D_i)}^{\tau_{\partial_j}^{'}\;+\tau_{\partial_j}^{'}}$
and $\theta_{D_i}$ (defined analogously as $\theta_{\partial_i}$)
is an automorphism of $\F_\infty\;\oplus\;\F_\infty\;\oplus\;S$.
Since $D_i$ is discontinuous, so is $\theta_{D_i}$ and so the
topologies $\tau_{\partial_j}^{'}\;+\tau_{\partial_j}^{'}$ and
$\tau_{D_i}^{'}$ (defined analogously as $\tau_{\partial_i}^{'}$)
are not equivalent on $\F_\infty\;\oplus\;\F_\infty\;\oplus\;S$.
It is easy to see that the topologies $\tau_{\partial_i}^{'}$ and
$\tau_{D_i}^{'''}$ (defined analogously as
$\tau_{\partial_i}^{'''}$) are equivalent on
$\F_\infty\;\oplus\;\F_\infty$ as $\tau_{\partial_i}^{'}\;\leq\;
\tau_{D_i}^{'''}$ on $\F_\infty\;\oplus\;\F_\infty$ and the open
mapping theorem for Fr\'{e}chet spaces. So, if the topologies
$\tau_{\partial_i}^{'}$ and $\tau_{\partial_j}^{'}$ are equivalent
on $\F_\infty\;\oplus\;\F_\infty$, then the topologies
$\tau_{\partial_j}^{'}\;+\tau_{\partial_j}^{'}$ and
$\tau_{D_i}^{'}$ are equivalent on
$\F_\infty\;\oplus\;\F_\infty\;\oplus\;S$, a contradiction. Hence
the Fr\'{e}chet algebra topologies $\tau_{\partial_i}^{'}$ are
mutually inequivalent.

Conversely, we show that $D_i$ is discontinuous on
$(\F_\infty\;\oplus\;\F_\infty,\;\tau_{\partial_j}^{'})$. For
this, let $((a_n, x_n))$ be a sequence in
$(\F_\infty\;\oplus\;\F_\infty,\;\tau_{\partial_j}^{'})$ such that
it converges to $(0, 0)$ in $\F_\infty\;\oplus\;\F_\infty$. So,
$q_{k, \partial_j}^{'}((a_n, x_n))\;\rightarrow\;0$. Now we see
that $D_i((a_n, x_n))\;=\;(0, \partial_i(a_n))$ does not converge
to $(0,\,0)$. Since $q_{k, \partial_j}^{'}((0,
\partial_i(a_n))\;=\;p_{k,\,j}^{'}(\partial_i(a_n))$, we show that $p_{k,\,j}^{'}(\partial_i(a_n))$ does not
converge to $0$. For this, suppose that for $i\,\neq\,j$, the
Fr\'{e}chet algebra topologies $\tau_{\partial_j}^{'}$ and
$\tau_{\partial_i}^{'}$ are not equivalent. Then, for all
$n\,\in\,\nn$ and $c\,>\,0$, there exists some $m\,\in\,\nn$ such
that
$$q_{m,\,\partial_j}^{'}((a,\,x))\,>\,c q_{n,\,\partial_i}^{'}((a,\,x))\;\;((a,\,x)\,\in\,\F_\infty\;\oplus\;\F_\infty).$$
Hence, for $(a,\,x)\,\in\,\F_\infty\;\oplus\;\F_\infty$,
$$q_{m,\,\partial_j}^{'}(D_i(a,\,x))=q_{m,\,\partial_j}^{'}((0,\,\partial_i(a)))=p_{m,\,j}^{'}(\partial_i(a))>c\,q_{n,\,\partial_i}^{'}((0,\,\partial_i(a))).$$
Since
$q_{n,\,\partial_i}^{'}((0,\,\partial_i(a)))=\,p_{n,\,i}^{'}(\partial_i(a))$
and $\partial_i$ is discontinuous on $(\F_\infty,\,\tau_i)$, we
see that $p_{k,\,j}^{'}(\partial_i(a_n))$ does not converge to $0$
whenever $a_n$ converges to $0$. So $D_i$ is discontinuous on
$(\F_\infty\;\oplus\;\F_\infty,\;\tau_{\partial_j}^{'})$. Hence,
the result follows. $\hfill \Box$\\

We remark that for each $i\,\in\,\zz^+$, the Fr\'{e}chet algebra
topologies $\tau_{\partial_i}$ and $\tau_{\partial_i}^{'}$ are
mutually inequivalent on the algebra
$\F_\infty\;\oplus\;\F_\infty$, because if they are equivalent,
then, by Corollary 3.7 and by the transitivity property of the
relation ``equivalence of topologies'', we see that for
$i\,\neq\,j$, the Fr\'{e}chet algebra topologies
$\tau_{\partial_i}$ and $\tau_{\partial_j}^{'}$ are equivalent on
the algebra $\F_\infty\;\oplus\;\F_\infty$, and so, for $i, \,
j\,\in\,\zz^+$ with $i\,\neq\,j$, the Fr\'{e}chet algebra
topologies $\tau_{\partial_i}^{'}$ and $\tau_{\partial_j}^{'}$ are
equivalent on the algebra $\F_\infty\;\oplus\;\F_\infty$, a
contradiction.

We give the following extension of Theorem ~\ref{Theorem 1_Loy} in
the Fr\'{e}chet case.
\begin{thm} \label{Thm. 2.1_Frechet} Let $A$ be a commutative Fr\'{e}chet
algebra, $(D_i)$ a sequence of non-zero derivations of $A$ into a
commutative Fr\'{e}chet $A$-module $M$. If $D_i$ vanishes on dense
subset of $A$ for each $i\,\in\,\zz^+$, then the algebra
$\overline{A}_{D_i}$ admits a Fr\'{e}chet algebra topology
$\tau_{D_i}^{'}$, generated by $(q_{k,\;D_i}^{'})$ (defined
analogously as $(q_{k,\;\partial_i}^{'})$), and is inequivalent to
the Fr\'{e}chet algebra topology, generated by $(q_{k,\,i}^{'})$
(defined analogously). $\hfill \Box$
\end{thm}
\section{Conclusions}
\noindent{\bf 4.1.} The above examples confirm the impression that
derivations can behave quite badly on Fr\'{e}chet algebras, which
impression is initially given simply by the fact that
$\partial_i\;=\;\partial/\partial X_i$ maps outside the radical,
on any formal power series algebra $\A_i[[X_i]]$, and so the
Singer-Wermer conjecture does not hold in this case. Also, it does
not hold for the derivation $f \mapsto f^{'}$ on the algebra
$\textrm{Hol}(U),\,U$ a domain in $\cc$. However this conjecture
(and its weaker version, too) does hold for the derivations $D_i$
on the Fr\'{e}chet algebra $(\F_\infty \oplus
\F_\infty,\,\tau_{\partial_i})$ as
$\textrm{Rad}(\F_\infty\,\oplus\,\F_\infty)\,=\,\textrm{Rad}(\F_\infty)\;\oplus\;\F_\infty$
and the image of $D_i$ is $0\,\oplus\,\F_\infty$; this also shows
that the stronger version (i.e., the Kaplansky's conjecture) does
not hold by this non-semisimple, non-Banach commutative
Fr\'{e}chet algebra. Moreover, its stronger version holds for
pro-$C^*$-algebras (and, in particular, for $C(X),\, X$ a
hemicompact $k$-space). So, it would be interesting to see which
classes of Fr\'{e}chet algebras satisfy these avatars of the
Singer-Wermer conjecture among all Fr\'{e}chet algebras (cf.
~\cite[\S 3, Question]{8}).

\noindent{\bf 4.2.} Moreover, it would also be interesting to
construct examples of Fr\'{e}chet algebras with uncountably many
inequivalent Fr\'{e}chet algebra topologies. The author believes
that the answer of this question is in the affirmative, but has
been unable to settle the question in general. For example, in the
Banach case, one may modify Feldman's example to have an example
of a Banach algebra with uncountably many inequivalent complete
algebra norms by using distinct (discontinuous) linear functionals
on a Banach algebra.\\

\noindent {\bf Acknowledgement.} A talk based on this paper was
delivered at the International Conference on Topological Algebras
and Their Applications (ICTAA 2021), held at Estonia, 7th to 9th
July, 2021. The author is grateful to the organizers of ICTAA 2021
for facilitating the talk.

\begin{tabbing}
\noindent Address: \= B-1/125, Swi-park Soc. (Madhuvrund),\\
\> Ghatlodia, Ahmedabad 380061, Gujarat, INDIA.
\end{tabbing}
\noindent E-mail: srpatel.math@gmail.com
\end{document}